# Multiple-rank Modification of Symmetric Eigenvalue Problem

HyungSeon Oh[*] and Zhe Hu

Electrical Engineering, University at Buffalo

**Abstract**

Rank-1 modifications in *k*-times (*k* > 1) often are performed to achieve rank-*k* modification. We propose a rank-*k* modification for enhancing computational efficiency. As the first step towards a rank-*k* modification, an algorithm to perform rank-2 modification is proposed and tested. The computation cost of our proposed algorithm is in $O(n^{1.5})$. We also propose a general rank-*k* update algorithm based upon the modified Sturm Theorem, and compare our results from those of the direct eigenvalue decomposition and of a perturbation method.

## 1. Introduction

Optimal power flow plays a key role in the operation and the planning studies for power systems. Due to its nonlinearity and non-convexity, a numerical solution is pursued using heuristic methods [1][2]. A most widely used method is a Lagrange relaxation, which relies on computationally expensive matrix factorizations. It is not uncommon that a matrix involves low-rank update. Therefore, it would be a viable option to update the factors instead of the expensive re-factorization process to evaluate them. An *LU* modification method is first introduced and many study results show its effectiveness [8][9]. While very relevant to power system studies because it can preserve the sparsity, the lack of numerically stability is an important issue. Stable update is achieved in updating Cholesky factorization [9-11], but its applicability is limited to positive semi-definite (PSD) matrices. The matrices associated with power systems are not, in general (PSD). Modification of eigenvalue/eigenvector pairs [12] would be a good candidate because they are numerically stable, and it can modify a large-scale matrix.

Consider a real symmetric matrix $A \in R^{n \times n}$ with known eigenvalue decomposition $A = Q \Lambda Q^T$, to which a symmetric perturbation is added. If the perturbation is a rank-one matrix, $\sigma v v^T$, the eigenvalue of new matrix $A + \sigma v v^T$ could be given by "secular equation" [4]

$$f(\lambda) = 1 + \sigma \sum_{j=1}^{N} \frac{\zeta_j^2}{\lambda_j - \lambda} = 0 \qquad (1)$$

where $\lambda_i$ is the *i*th eigenvalue of original matrix *A*, and $\zeta_i$ is the *i*th element in vector $z = Q^T v$. If the $\zeta_i$'s are all nonzero and $\lambda_i$ are distinct, then this equation has n solutions. On interval $[\lambda_i, \lambda_{i+1}]$, function f is monotonic. These properties lead to several efficient and stable methods to solve the secular equation [4][5].

However, it is difficult to find the rank-1 modification matrix during the heuristic algorithm. In most situations, the perturbation is not a rank-1 matrix. A reduction of the rank in the perturbation matrix is possible by utilizing eigenvalue decomposition after selecting a subset of eigen-pairs. The disadvantages of this approach are: 1. computationally expensive eigenvalue decomposition; 2. poor accuracy when a subset of eigenvalue pairs are included. It will be beneficial to modify high-rank eigenvalue directly. Eigenvalue decomposition is performed for a symmetric matrix, and the perturbation is made in preserving the symmetry. In an iterative method, binding constraint sets changes, which involves an update such as $e_j^T a + a^T e_j$. Such a change is easily recognized without any further analysis to identify the perturbation vector as in rank-1 update process. This is a main motivation for rank-2 update.

In this paper, we propose a new, direct rank-$k$ modification where $k \geq 2$. In Section 2, we present the theoretical background for the multiple-rank modification for rank-$k$ updates. Section 3 lists the numerical results, and Section 4 outlines conclusions and future works.

## 2. Theory of Eigenvalue Updates

### 2.1. Secular equation of multiple-rank modification

Multiple-rank modification theory was carried out in [3]. Considering the eigenvalue decomposition of a real symmetric matrix $A + KK^T$, where $K \in R^{n \times k}$, $k \ll n$. The eigenvalue decomposition of matrix $A$ is known. It is pointed out that the solution of multiple-rank modification problem is equivalent to roots of

$$f(\lambda) = \det\left[I_r - K^T(\lambda I - A)^{-1}K\right] = \det\left[I_r - U^T(\lambda I - \Lambda)^{-1}U\right] = 1 + \sum_{k=1}^{r}(-1)^k \sum_{\substack{1 \leq r_1 \leq \cdots \leq r_k \leq r \\ 1 \leq j_1 \leq \cdots \leq j_k \leq n}} \frac{\left[\det(U_{j_1 \cdots j_k; r_1 \cdots r_k})\right]^2}{(\lambda - \lambda_{j_1}) \cdots (\lambda - \lambda_{j_k})} \quad (1)$$

where $U = Q^T K$, $U_{j_1 \cdots j_k; r_1 \cdots r_k} = \begin{pmatrix} e_{j_1} & \cdots & e_{j_k} \end{pmatrix}^T U \begin{pmatrix} e_{j_1} & \cdots & e_{j_k} \end{pmatrix}$, and $e_j$ is the $j^{th}$ column vector from an identity matrix. Golub point out that solving $f(\lambda) = 0$ is numerical difficult [3]. While calculating the value of determinants is expensive, rearranging it is possible in a following way:

$$\sum_{k=1}^{r}(-1)^k \sum_{\substack{1 \leq r_1 \leq \cdots \leq r_k \leq r \\ 1 \leq j_1 \leq \cdots \leq j_k \leq n}} \frac{\left[\det(U_{j_1 \cdots j_k; r_1 \cdots r_k})\right]^2}{(\lambda - \lambda_{j_1}) \cdots (\lambda - \lambda_{j_k})} = \sum_{k=1}^{r}(-1)^k \sum_{\substack{1 \leq r_1 \leq \cdots \leq r_k \leq r \\ 1 \leq j_1 \leq \cdots \leq j_k \leq n}} \left[\det(U_{j_1 \cdots j_k; r_1 \cdots r_k})\right]^2 \sum_{i \in \{j_1 \cdots j_k\}} \left[\frac{1}{\prod_{\substack{s \neq i \\ s \in \{j_1 \cdots j_k\}}}(\lambda_s - \lambda_i)}\right] \frac{1}{\lambda_i - \lambda}$$

$$= \sum_{i=1}^{N} \frac{1}{\lambda_i - \lambda} \left\{ \sum_{k=1}^{r} \sum_{\substack{1 \leq r_1 \leq \cdots \leq r_k \leq r \\ 1 \leq j_1 \leq \cdots \leq j_k \leq n}} \frac{\left[\det(U_{j_1 \cdots j_k; r_1 \cdots r_k})\right]^2}{\prod_{\substack{s \neq i \\ s \in \{j_1 \cdots j_k\}}}(\lambda_s - \lambda_i)} \right\} = \sum_{i=1}^{N} \frac{\alpha_i^{\{0\}}}{\lambda_i - \lambda}$$

where $D_i = diag\left(\frac{1}{\lambda_1 - \lambda_i}, \cdots, \frac{1}{\lambda_{i-1} - \lambda_i}, 0, \frac{1}{\lambda_{i+1} - \lambda_i}, \cdots, \frac{1}{\lambda_n - \lambda_i}\right)$ and $\alpha_i^{\{0\}} = \det\begin{bmatrix} 0 & e_i^T U \\ U^T e_i & I_r + U^T D_i U \end{bmatrix}$ (2)

Therefore, similar to rank-1 modification, $f$ also can have a secular equation,

$$f(\lambda) = 1 + \sum_{j=1}^{N} \frac{\alpha_i^{\{0\}}}{\lambda_j - \lambda} \quad (3)$$

This secular equation is called "multiple-rank secular equation". The matrix in determinant is $(k+1) \times (k+1)$ matrix. When $k \ll N$, formulating $\mu_i$ will not increase the computation cost significantly. When $k$ increases, the determinants calculation involves a heavy computation, which makes this approach inefficient.

Another advantage is that this method does not require an orthogonal perturbation matrix $V$. Therefore, an additional eigenvalue decomposition for the perturbation matrix is not necessary. However, solving multiple-rank secular equations is more computationally demanding than that in the rank-1 update process. In rank-1 modification, $\mu_i$ are always positive, and $f$ is monotonically increasing on each sub-interval $(\lambda_i, \lambda_{i+1})$. Therefor only one root exists in the range of $[\lambda_i, \lambda_{i+1}]$. In updating multiple-rank, $\mu_i$ are not necessarily positive, so that $f$ is not monotonic on $(\lambda_i, \lambda_{i+1})$. There would be at least 0 and at most $r$ solutions. Although multiple-rank modification problems have similar secular equation to rank-1, most methods used to solve rank-1 problems in [3] cannot be directly used for multiple-rank problems. To solve the multiple-rank secular equations, we formulate two sub-problems: the first sub-problem intends to locate eigenvalues in checking how many roots are on each sub-interval; the second one finds the value for these roots on each sub-interval. The details are listed in Section 3. We may not need to update all the eigenvalues, if only the eigenvalues in some specific interval are of interest, then the first sub-problem is enough.

## 2.2. Location of eigenvalues

In the process for the rank-2 modification, eigenvalues are located according to the signs of $\mu_i$. For the case of rank-$k$, they become more complicated to locate. In general, the one-column-one-row update approximates the modification while it preserves the sparsity and the symmetry [5]:

Step 1. $A = A_0 + ab^T + ba^T = Q(\Lambda + xy^T + yx^T)Q^T$

Step 2. $\Lambda + xy^T + yx^T = \Lambda + (x \quad y)\begin{bmatrix} 0 & 1 \\ 1 & 0 \end{bmatrix}(x \quad y)^T = \Lambda + (x \quad y)S\begin{bmatrix} 1 & 0 \\ 0 & -1 \end{bmatrix}S^T(x \quad y)^T$

where $S$ is the Schur compliment. Let $(u_1 \quad u_2) = (x \quad y)S$, which yields $\Lambda + xy^T + yx^T = \Lambda + u_1 u_1^T - u_2 u_2^T$. Therefore, if we can find out the dominant column-row pair, it is acceptable to update eigenvalues using several rank-2 modification.

Considering the secular equation of rank-2 modification $A + \sigma_1 v_1 v_1^T + \sigma_2 v_2 v_2^T$, with $\sigma_1 > 0$ and $\sigma_2 > 0$. In the case, the values of all the eigenvalues increase, termed double-right shift. $R$ is divided into $N+1$ sub-intervals by original eigenvalues $\{\lambda_1, \dots, \lambda_N\}$. When the rank of the perturbation matrix is 2, there are 0,1 or 2 roots within each sub-interval. We define a location vector $L \in Z^{(n+1)\times 1}$, where each element represents the number of roots within the corresponding sub-interval. We can approximate locations of new eigenvalues using Courant-Weyl principle [6],

$$\lambda_k \leq \lambda_k^{\max}; \quad 1 \leq k \leq n \text{ and } \lambda_k^{\max} \leq \lambda_{k+2}; \quad 1 \leq k \leq N-2$$

Also, from the signs of $\mu_i$, parity of $L$ could be determined,

*If $\mu_i \cdot \mu_{i+1} > 0$,* there is one eigenvalue in $(\lambda_i, \lambda_{i+1})$; otherwise, there are either zero or two eigenvalues in this interval $(\lambda_i, \lambda_{i+1})$. These three conditions are already sufficient to locate new eigenvalues. Figure 1 illustrates a pseudo-code of Algorithm 1a for the location vector.

---
Algorithm 1a. Locating eigenvalues (double-right shift)

$M = 0; L_1 = 0$

for $i = 1:n$
    if $\mu_i \cdot \mu_{i+1} > 0$
        $L_{i+1} = 1; M = M+1;$
    else if $M+2 > i$
        $L_{i+1} = 0;$
    else
        $L_{i+1} = 2; M = M+2;$
    end
end

Figure 1. Pseudo-code of formulating locating vector when the matrix is under double-right shift update.

Note that the condition of Algorithm 1a is $\sigma_1 > 0$ and $\sigma_2 > 0$. The condition can be generalized either $\sigma_1 < 0$ and $\sigma_2 < 0$ ; or $\sigma_1 < 0$ and $\sigma_2 > 0$. For the first case, instead of searching for a root from the first subinterval, we can start from the last subinterval. And for situation (c), sign of $\mu_1$ can determine number of eigenvalues on first interval. Figure 2 shows the pseudo-code for these to condition.

---
**Algorithm 1b. Locating eigenvalues (double-left shift)**

$M = 0; L_{n+1} = 0$

for $i$ = 1:n
    if $\mu_{n-i} \cdot \mu_{n-i+1} > 0$
        $L_{n-i+1} = 1; M = M + 1;$
    else if $M + 2 > n - i + 1$
        $L_{n-i+1} = 0;$
    else
        $L_{n-i+1} = 2; M = M + 2;$
    end
end

---
**Algorithm 1c. Locating eigenvalues (one-right-one-left shift)**

if $\mu_1 > 0$
  $M = 0; L_1 = 0$
else
  $M = 1; L_1 = 1$
end

for $i$ = 1:n
    if $\mu_i \cdot \mu_{i+1} > 0$
        $L_{i+1} = 1; M = M + 1;$
    else if $M + 2 > i$
        $L_{i+1} = 0;$
    else
        $L_{i+1} = 2; M = M + 2;$
    end
end

---

Figure 2. Pseudo-code of generating location vector when matrix is under double-left shift update and one-left-one-right shift update.

For general rank-$k$ modification sign of $\mu_i$ is not sufficient to calculate the location vector. We propose a modified Sturm Chain method as follows.

The original Sturm chain or Sturm sequence [7] is a finite sequence of polynomials $p_0(\lambda), p_1(\lambda), \cdots, p_m(\lambda)$ of decreasing degree with these following properties:

- $p_0(\lambda)$ is square free (no square factors, i.e., no repeated roots);
- $p_0(\lambda) = 0$, then $\text{sign}(p'(\lambda)) = \text{sign}(p_1(\lambda))$
- If $p_i(\lambda) = 0, 0 < i < m$, $\text{sign}(p_{i-1}(\lambda)) = -\text{sign}(p_{i+1}(\lambda))$
- $p_m(\lambda)$ does not change its sign.

Then by observing signs of $p_0(a), p_1(a) \ldots p_m(a)$ and $p_0(b), p_1(b) \ldots p_m(b)$, number of the roots of $p_0(\lambda)$'s within (a, b) could be determined. Although original Strum Theorem is for polynomial function, we found that it can also be implemented on secular function. For our secular equation, define that,

$$f_0(\lambda) = 1 - \sum_{j=1}^{N} \frac{\alpha_i^{\{0\}}}{\lambda - \lambda_j} \tag{4}$$

Note that $d_i$ is not repeated. For $f_1(\lambda)$, in order to satisfy the second condition, we need to take the derivative of $f_0(\lambda)$. Direct derivation of $f_0(\lambda)$ will lead to second order terms. Instead, take the derivative of following function,

$$p_0(\lambda) = \pi_0(\lambda)\left(1 - \sum_{j=1}^{N} \frac{\alpha_i^{\{0\}}}{\lambda - \lambda_j}\right) \text{ where } \pi_m(\lambda) = \prod_{j=1}^{N-m}(\lambda - \lambda_j) \tag{5}$$

Then we have,

$$\begin{aligned}
p_1(\lambda) &= \frac{dp_0(\lambda)}{d\lambda} = \frac{df_0(\lambda)}{d\lambda}\pi_0(\lambda) + f_0(\lambda)\left(\sum_{j=1}^{N}\frac{1}{\lambda - \lambda_j}\right)\pi_0(\lambda) \\
&= \left\{\sum_{j=1}^{N}\frac{\alpha_i^{\{0\}}}{\lambda - \lambda_j}\frac{1}{\lambda - \lambda_j} + \sum_{j=1}^{N}\frac{1}{\lambda - \lambda_j} - \sum_{j=1}^{N}\sum_{i=1}^{N}\frac{1}{\lambda - \lambda_i}\frac{\alpha_i^{\{0\}}}{\lambda - \lambda_j}\right\}\pi_0(\lambda) \\
&= \left\{\sum_{j=1}^{N}\frac{\alpha_i^{\{0\}}}{\lambda - \lambda_j}\frac{1}{\lambda - \lambda_j} + \sum_{j=1}^{N}\frac{1}{\lambda - \lambda_j} - \sum_{j=1}^{N}\sum_{i=1}^{N}\frac{1}{\lambda - \lambda_i}\frac{\alpha_i^{\{0\}} + \alpha_j^{\{0\}}}{\lambda_j - \lambda_i}\right\}\pi_0(\lambda) = \left[\sum_{j=1}^{N}\frac{\mu_j^{\{1\}}}{\lambda - \lambda_j}\right]\pi_0(\lambda) \\
&= \left[c_1 - \sum_{j=1}^{N-1}\frac{\alpha_j^{\{1\}}}{\lambda - \lambda_j}\right]\pi_1(\lambda) \text{ where } \mu_j^{\{1\}} = 1 - \sum_{\substack{j=1 \\ i \neq j}}^{N}\frac{\alpha_i^{\{0\}} + \alpha_j^{\{0\}}}{\lambda_j - \lambda_i}, \ c_1 = \sum_{j=1}^{N}\mu_j^{\{1\}}, \text{ and } \alpha_j^{\{1\}} = (\lambda_N - \lambda_j)\mu_j^{\{1\}}
\end{aligned} \tag{6}$$

Here, $\alpha_i^{\{1\}}$ is not related with $\lambda$. The Sturm chain theorem finds the number of solutions inside a given range to satisfy $p_0(\lambda) = 0$. However, the direct application of the theorem to $p_0(\lambda)$ is not practical because the computation of the coefficients is numerically unstable. Inspired by the fact that only the change in the signs of $p_i(\lambda)$ is important, which is the product between $\pi_i(\lambda)$ and $f_i(\lambda)$, we propose a modified Sturm series in a secular form: The following elements $p_2(\lambda) \ldots p_N(\lambda)$ are generated by long division used for polynomial function that can be also expanded to secular function.

---

**Secular Long Division**

For secular functions with same eigenvalues $\lambda$,

$$p_m(x) = \left(c_m - \sum_{j=1}^{N-m}\frac{\alpha_j^{\{m\}}}{x - d_j}\right)\prod_{k=1}^{N-m}(x - d_k), \ p_{m-1}(x) = \left(c_{m-1} - \sum_{j=1}^{N-m+1}\frac{\alpha_j^{\{m-1\}}}{x - d_j}\right)\prod_{k=1}^{N-m+1}(x - d_k), \text{ and}$$

$$p_{m-2}(x) = \left(c_{m-2} - \sum_{j=1}^{N-m+2}\frac{\alpha_j^{\{m-2\}}}{x - d_j}\right)\prod_{k=1}^{N-m+2}(x - d_k), \text{ there exists secular long division}$$

$f_{m-2}(x)\pi_{m-2}(x) - f_{m-1}(x)A_m(x - B_m)\pi_{m-1}(x) + f_m(x)\pi_m(x) = 0$ where

$A_m = \frac{c_{m-2}}{c_{m-1}}$ and $B_m = -d_{N-m+2} - \frac{1}{c_{m-2}}\sum_{j=1}^{N-m+2}\alpha_j^{\{m-2\}} + \frac{1}{c_{m-1}}\sum_{j=1}^{N-m+1}\alpha_j^{\{m-1\}}$. With a given $A_m$

and $B_m$, one finds $p_m(x) = f_m(x)\pi_m(x) = -f_{m-2}(x)\pi_{m-2}(x) + f_{m-1}(x)A_m(x - B_m)\pi_{m-1}(x)$.

---

In the Sturm chain theorem, the number of changes in sign of the polynomial functions $p_0(\lambda) \ldots p_N(\lambda)$ equals the number of solutions inside a region. In $[-\infty, \infty]$, the signs of $f_0(\lambda) \ldots f_N(\lambda)$ are same as those of $c_0 \ldots c_N$; the signs of $\pi_0(\infty) \ldots \pi_N(\infty)$ are all positive; and those of $\pi_0(-\infty) \ldots \pi_N(-\infty)$ are all negative. As a result, the number of solutions that satisfy $p_0(\lambda) = 0$ is number of positives in $c_0 \ldots c_N$ where $c_0 = 1$. In the computation process of $c_0 \ldots c_N$, the floating number can result in a change in the sign of $c_0 \ldots c_N$. If one terminates the Secular long division when it hits the non-positive $c$ for the first time, the partial Sturm series yields the lower bound of the solutions between $[-\infty, \infty]$. Furthermore, one can break down the regions where at least a solution exists, which yields a solution as will be described in Section 2.3. Suppose, $(N_s^0 + 1)$ $c$'s are positive at the first Sturm series, and accordingly $N_s^0$ solutions $\xi_s$ are identified. Then a new polynomial function $p_0^1(\lambda)$ is defined as follows:

$$p_0(\lambda) = \left[1 - \sum_{j=1}^{N} \frac{\alpha_j^{\{0\}}}{\lambda - \lambda_j}\right] \prod_{j=1}^{N}(\lambda - \lambda_j) = p_0^1(\lambda) \prod_{s=1}^{N_s^0}(\lambda - \xi_s) \rightarrow p_0^1(\lambda) = \left[1 - \sum_{j=1}^{N} \frac{\alpha_j^{\{0\}}}{\lambda - \lambda_j}\right] \frac{\prod_{j=1}^{N}(\lambda - \lambda_j)}{\prod_{s=1}^{N_s^0}(\lambda - \xi_s)} \quad (7)$$

By multiplying the first term after the square parenthesis, (7) becomes:

$$p_0^1(\lambda) = \left[1 - \sum_{j=1}^{N} \frac{\alpha_j^{\{0\}}}{\lambda - \lambda_j}\right] \frac{\lambda - \lambda_1}{\lambda - \xi_1} \frac{\prod_{j=2}^{N}(\lambda - \lambda_j)}{\prod_{s=2}^{N_s^0}(\lambda - \xi_s)} = \left[1 - \sum_{j=1}^{N}\left(\frac{d_j - d_1}{d_j - \xi_1}\right)\frac{\alpha_j^{\{0\}}}{x - d_j} - \frac{d_1 - \xi_1}{x - \xi_1}\left(1 - \sum_{j=1}^{N}\frac{\alpha_j^{\{0\}}}{\xi_1 - d_j}\right)\right] \frac{\prod_{j=2}^{N}(\lambda - \lambda_j)}{\prod_{s=2}^{N_s^0}(\lambda - \xi_s)} \quad (8)$$

Since $\xi_1$ is a solution to $p_0(\lambda) = 0$ (i.e., $p_0(\xi_1) = 0$) but $\pi_0(\xi_1) \neq 0$, $1 - \sum_{j=1}^{N} \frac{\alpha_j^{\{0\}}}{\xi_1 - d_j}$ vanishes. After multiplying all $N_s^0$ terms, $p_0^1(\lambda)$ becomes:

$$p_0^1(\lambda) = \left[1 - \sum_{j=1}^{N} \frac{\beta_j^{\{0\}}}{\lambda - \lambda_j}\right] \prod_{j=N_s^0+1}^{N}(\lambda - \lambda_j) \text{ where } \beta_j^{\{0\}} = \left[\prod_{p=1}^{N_s^0} \frac{d_j - d_p}{d_j - \xi_p}\right] \alpha_j^{\{0\}} \quad (9)$$

The formula for $p_0^1(\lambda)$ resembles that for $p_0(\lambda)$, which allows the Sturm series expansion with less terms by $N_s^0$. Since the first $c$ in $p_0^1(\lambda)$ is positive, it is guaranteed to have at least one solution using the new Sturm series. This process will continue until all $N$ solutions are identified.

## 2.3. Algorithm to find solutions for the secular equation

Different from rank-1 modification, the secular equation of rank-2 is not monotonously increasing in each interval. Therefore, most of the algorithms for rank-1 update may not converge. Divide and Conquer (D&C) method is efficient to address this problem. As long as location vector is formulated, D&C is fast and parallelizable. Fig 3 is an example of rank-2 secular equation. By using algorithm in Section 2.2, we can get the location vector.

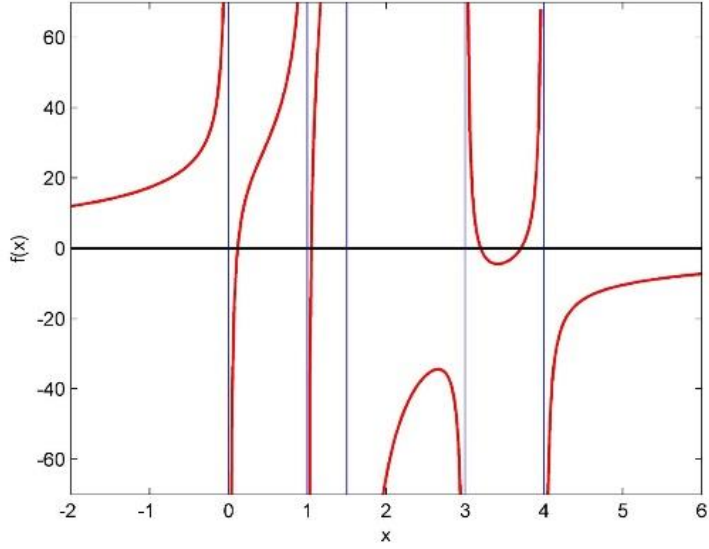

Figure 3, Graph of rank-2 secular function

With location vector, D&C method can be utilized.

**Divide-and-Conquer zero finder**
1. Choose $m$ equidistant points between upper bound $\overline{x}$ and lower bound $\underline{x}$, where a solution exists.
2. Check the signs of $f(x)$ at each point. Decide bounds of solutions.
3. Update the $\overline{x}$ or $\underline{x}$.
4. If $|\overline{x} - \underline{x}| < tol$, terminate.

## 2.4. Algorithm to update eigenvectors

Once an eigenvalue is computed, the corresponding eigenvector $v$ can be update using the relationship between an eigenvalue and an eigenvector, i.e., $(\Lambda + LL^T)u = \lambda u$ where $L = Q^T K$ and $u = Q^T v$. Solving the relationship for $v$ yields $u = (\Lambda - \lambda I)^{-1} L\alpha$ where $\alpha = -L^T u, \alpha \in R^k$ that implies the eigenvector is in the column space spanned by the rank-$k$ update $K$. α is in the null space of $\{I + L^T(\Lambda - \lambda I)^{-1} L\}^{k \times k}$, and let $α_0$ denote the null space vector. The row rank update $k$ where $k \ll n$ males finding the null space not computationally expensive. Since the eigenvector is a unit vector, $v$ becomes.

$$v = \frac{Q(\Lambda - \lambda I)^{-1} Q^T K\alpha_0}{\left\| Q(\Lambda - \lambda I)^{-1} Q^T K\alpha_0 \right\|_2} \tag{10}$$

## 3. Numerical Result

Figure 4 illustrates the comparison of the computation time of rank-1 update, rank-2 update, and original eigenvalue decomposition. The testing results shows the complexities $p$ of the computation are in $O(n^p)$:

- Two times rank-1 update: $p = 1.81$

- Rank-2 update: $p = 1.45$
- Direct EVD: $p = 2.88$

The computation cost of rank-2 modification is significantly better than that of rank-1 modification. Parallel computation will further enhance the computational efficiency.

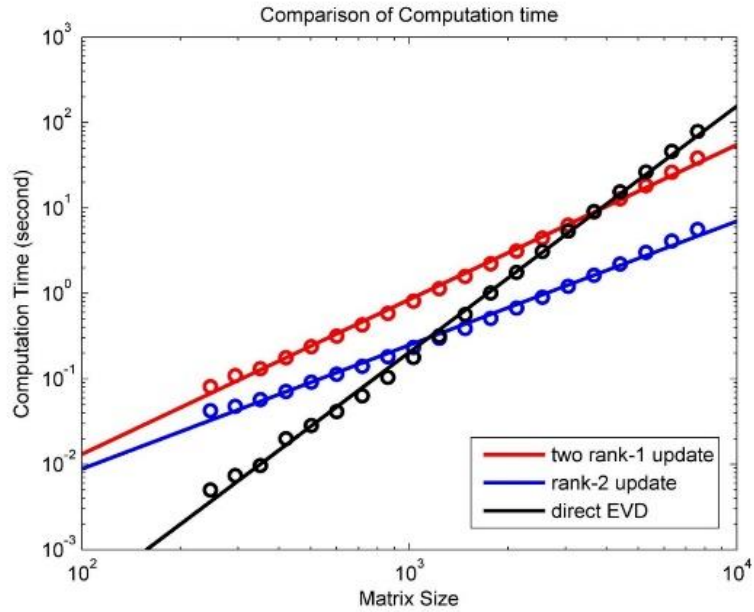

Figure 4, Computation time of two times rank-1 update, rank-2 update, and original eigenvalue decomposition

In Figure 5, the red line is the computation time for a randomly generated rank-3 update, and blue line is that of original eigenvalue decomposition. Updating the eigenvalues by secular Sturm Chain method is in $O(n^{1.95})$, which is more expensive than rank-2 update, but still much more efficient than original eigenvalue decomposition.

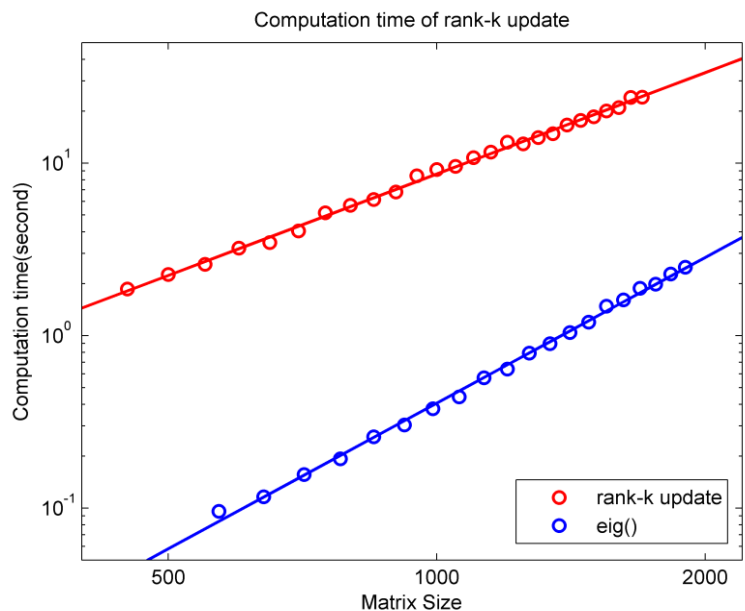

Figure 5. computation time of formulating locating vector by secular Sturm Chain.

When the perturbed matrix $K$ is small, i.e., norm($K$) << 1, the eigenvalue and eigenvector pairs can be estimated from:

$$\lambda_j = \lambda_j^0 + \left(K^T v_j^0\right)^T \left(K^T v_j^0\right) \text{ and } v_j = v_j^0 + \sum_{\substack{i=1 \\ i \neq j}}^{N} \frac{\left(K^T v_j^0\right)^T \left(K^T v_i^0\right)}{\lambda_i^0 - \lambda_j^0} v_j^0 \qquad (11)$$

Figure 6 illustrates the comparisons in the results among direct eigenvalue decomposition, proposed method, and the perturbation method when the norm($K$) is 0.01 and $N = 100$. When the matrix $K$ is very small the perturbation method provides a reasonably close estimate of the eigenvalue and eigenvector pair.

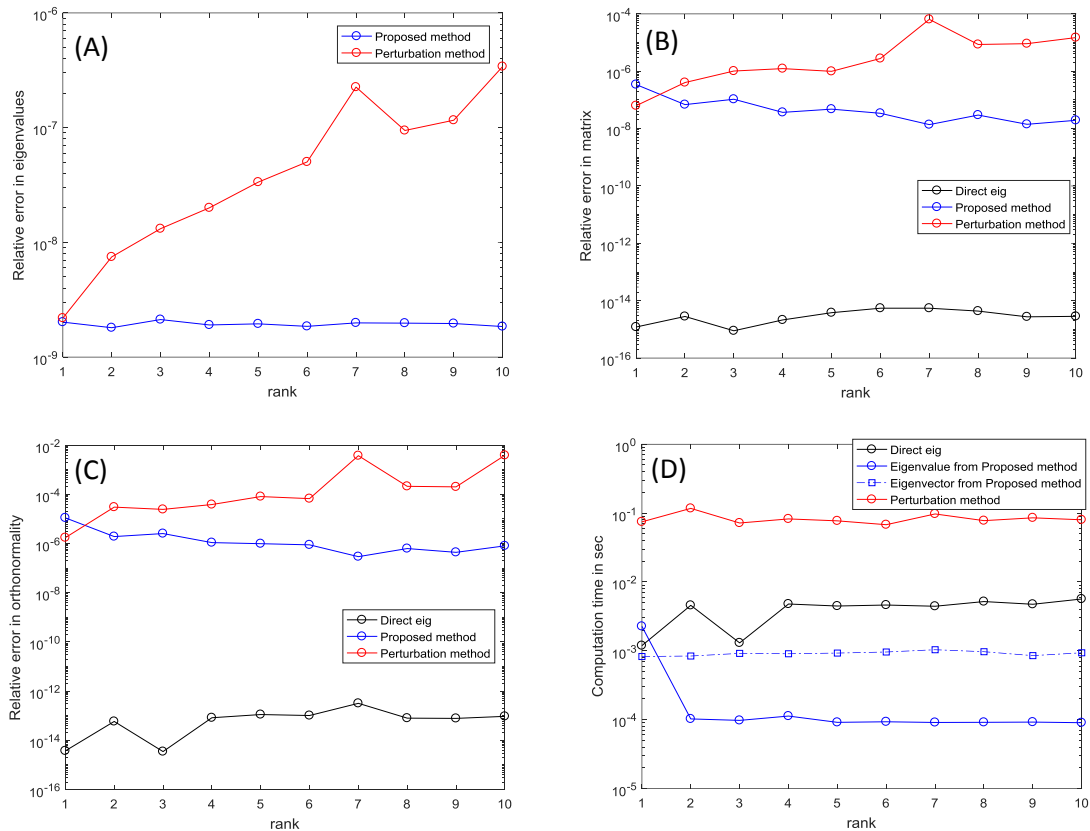

Figure 6. Comparison of the results in terms of the rank of $K$ from direct eigenvalue decomposition, proposed method, and the perturbation method outlined in (11) when norm($K$) = 0.01. (A) 2-norm error in the eigenvalues relative to the direct decomposition, (B) 2-norm error of the matrix before and after the decomposition, (C) 2-norm error of the products of eigenvector matrix and its transpose to check the orthonormality, and (D) Computation time of direct eigenvalue decomposition, of the eigenvalue update and of the eigenvector update of the proposed method, and of the perturbation method.

However, as the norm increases, the perturbation method starts to fail since its assumption to ignore higher order terms is not valid. Figure 7 shows the comparison with the norm($K$) is 0.3 and $N = 100$. While the performance of the proposed method does not change with the norm of $K$, the eigenvalues, the orthonormality, and the product between eigenvalue and eigenvector pairs all deviate from the desired quantities.

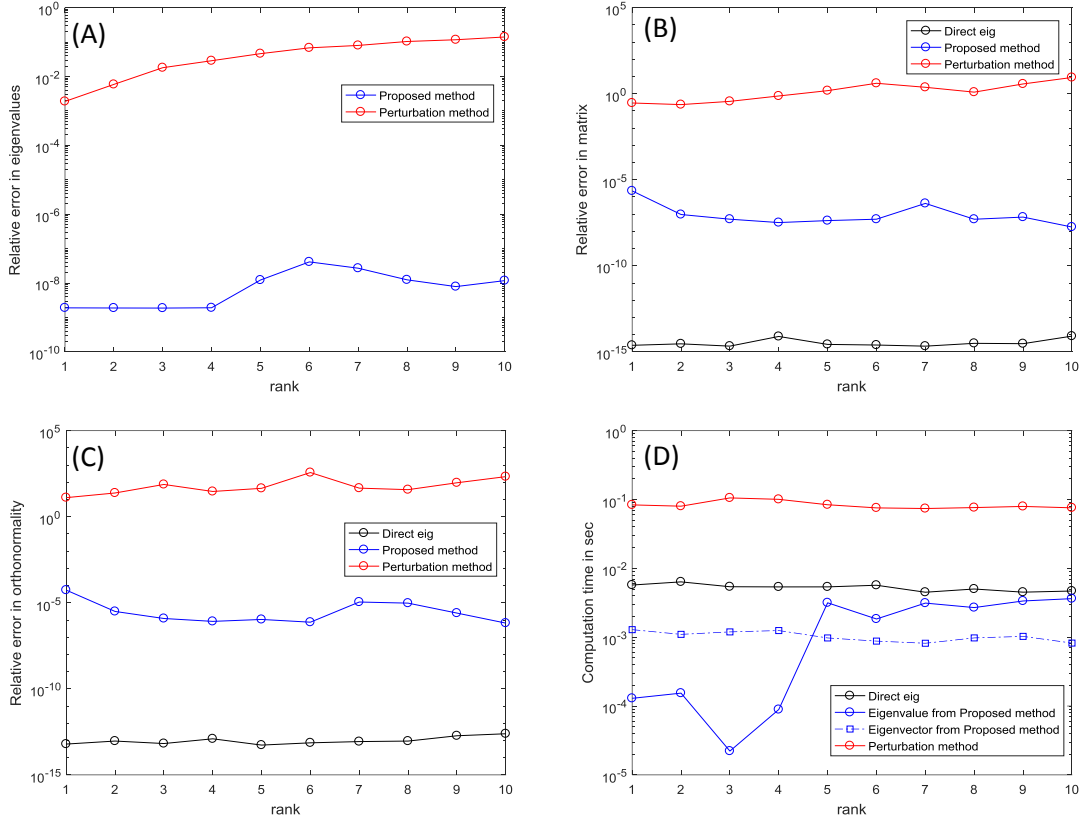

Figure 7. Comparison of the results in terms of the rank of *K* from direct eigenvalue decomposition, proposed method, and the perturbation method outlined in (11) when norm(*K*) = 0.3. (A) 2-norm error in the eigenvalues relative to the direct decomposition, (B) 2-norm error of the matrix before and after the decomposition, (C) 2-norm error of the products of eigenvector matrix and its transpose to check the orthonormality, and (D) Computation time of direct eigenvalue decomposition, of the eigenvalue update and of the eigenvector update of the proposed method, and of the perturbation method.

While the proposed method is efficient when the rank of the matrix *K* is small (approximately 10% of *N*), the computation time increases as the rank of *K* increases. Even if the computation time of the proposed method becomes larger than that of direct eigenvalue decomposition, the proposed method is still efficient to compute a subset of eigenvectors. For example, a null space is useful in many engineering fields and the identification of an entire eigenvector space may not be necessary. In such a circumstance, the eigenvalues are all zeros and the computation time to find eigenvectors are still much smaller than that of direct eigenvalue decomposition (compare the dotted blue line and the black solid line in (D) from Figures 6 and 7).

## 4. Conclusion

In this paper, our new algorithm for a rank-*k* modification of eigenvalue decomposition is presented. Computation performance is significantly improved in comparison to the rank-1 modification methods. We proposed a method relying on the location vector and tested on several systems. The test results show our proposed method is efficient.

**Acknowledgement:**

This work is supported by the Depart of Energy through the CERTS (Consortium for Electric Reliability Technology Solutions) program, PO# 7069673.